\documentclass{article}
\usepackage[utf8]{inputenc}
\usepackage{amsmath,amsfonts,amssymb,amsthm,mathtools}
\usepackage{todonotes}
\usepackage{hyperref}
\usepackage{csquotes}
\usepackage{enumerate}

\usepackage{geometry}
\newcommand{\LCP}{\operatorname{LCP}}
\newcommand{\define}{\coloneqq}
\newcommand{\bbR}{\mathbb{R}}
\newcommand{\bbRn}{\mathbb{R}^n}

\newcommand{\SnMi}{\mathbb{S}^{n-1}}
\newcommand{\sign}{\operatorname{sign}}
\newcommand{\bF}{\bar F}

\newcommand{\MnR}{\operatorname{M}_n(\mathbb R)}
\DeclareMathOperator{\cond}{\kappa^a}
\DeclareMathOperator{\dist}{\mathrm{dist}}

\newcommand{\Signs}{\mathcal{S}_n}
\newcommand{\spec}{\operatorname{Spec}}
\DeclareMathOperator{\specSR}{\operatorname{Spec}^{\rho}}
\DeclareMathOperator{\specA}{\operatorname{Spec}^a}
\DeclareMathOperator{\specLCP}{\operatorname{Spec}^a_{\operatorname{L}}}
\DeclareMathOperator{\countA}{\operatorname{c}^a}
\DeclareMathOperator{\countLCP}{\operatorname{c}^a_L}

\DeclareMathOperator{\img}{\operatorname{im}}

\newcommand{\eproof}{\hfill\qed}

\newtheorem{lemma}{Lemma}[section]
\newtheorem{corollary}[lemma]{Corollary}
\newtheorem{theorem}[lemma]{Theorem}
\newtheorem{proposition}[lemma]{Proposition}
\theoremstyle{definition}
\newtheorem{definition}[lemma]{Definition}
\theoremstyle{remark}
\newtheorem{remark}[lemma]{Remark}
\newtheorem{example}[lemma]{Example}

\title{Generalized Perron Roots\\and Solvability of the Absolute Value Equation}
\author{Manuel Radons\\
Technische Universität Berlin\\
Chair of Discrete Math./Geometry\\
Berlin, GERMANY\\
{\tt radons@math.tu-berlin.de}
\and
Josué Tonelli-Cueto\\
Inria Paris \& IMJ-PRG\\
OURAGAN team\\
Paris, FRANCE\\
{\tt josue.tonelli.cueto@bizkaia.eu}
}
\date{ }

\usepackage[
backend=biber,
style=numeric,
sorting=nty,
maxbibnames=99
]{biblatex}
\AtEveryBibitem{\clearfield{issn}}
\AtEveryBibitem{\clearfield{isbn}}
\AtEveryBibitem{\clearfield{url}}
\AtEveryBibitem{\ifentrytype{book}{\clearfield{pages}}{}}

\begin{document}
	
\maketitle
	
\begin{abstract}
Let $A$ be a $n\times n$ real matrix. The piecewise linear equation system $z-A\vert z\vert =b$ is called an absolute value equation (AVE). It is well-known to be equivalent to the linear complementarity problem. Unique solvability of the AVE is known to be characterized in terms of a generalized Perron root called the sign-real spectral radius of $A$. For mere, possibly non-unique, solvability no such characterization exists. We narrow this gap in the theory. That is, we define the concept of the aligned spectrum of $A$ and prove, under some mild genericity assumptions on $A$, that the mapping degree of the piecewise linear function $F_A:\mathbb{R}^n\to\mathbb{R}^n\,, z\mapsto z-A\lvert z\rvert$ is congruent to $(k+1)\mod 2$, where $k$ is the number of aligned values of $A$ which are larger than $1$. We also derive an exact---but more technical---formula for the degree of $F_A$ in terms of the aligned spectrum. Finally, we derive the analogous quantities and results for the LCP.
\end{abstract}

\section{Introduction}

	The \emph{linear complementarity problem} $\operatorname{LCP}(q,M)$, where $q\in\bbRn$, and $M\in\MnR$, the space of $n\times n$ real matrices, is to determine $v,w\in\bbRn_{\geq 0}$ with $v^Tw=0$ so that
	\begin{equation}\label{eq:LCP-standard}
	    v\ =\ Mw + q\,.
	\end{equation}
	It provides a common framework for numerous optimization tasks in economics, engineering and computer science. Classical problems that can be reduced to solving an LCP include bimatrix games, and linear and quadratic programs \cite{cottle1992lcp}. 
	Recent applications are the correct formulation of numerical models for free-surface hydrodynamics \cite{brugnano2008iterative}, $L_1$ regularization in reinforcement learning \cite{NIPS2010_81dc9bdb}, and the massively parallel implementation of collision detection on CUDA GPUs \cite[Chap. 33]{nguyen2007gpu}.
	
	It is well known that an $\operatorname{LCP}(q,M)$ is uniquely solvable for arbitrary $q$ if and only if $M$ is a $P$-matrix, that is, a matrix whose principal minors are all positive.
	Checking whether a matrix is a $P$-matrix is \texttt{co-NP}-complete~\cite{coxson1994p}. 
	As a consequence, there exists a rich body of literature about sufficient criteria for unique solvability, e.g., in terms of the matrix structure \cite[Thm. 3.2]{rump2003p} (cf.~\cite[Chap. 3]{cottle1992lcp}, \cite{zamani2022condition}), and in terms of norm constraints \cite[Thm. 2.15]{rump1997theorems} (cf.~\cite[Thm. 3.1]{rump2003p}). 
	
	If $\operatorname{LCP}(q,M)$ is solvable---possibly non-uniquely---for arbitrary $q$, then we call $M$ a \emph{$Q$-matrix}. 
	Unlike $P$-matrices, $Q$-matrices have quite involved characterizations~\cite{deloeramorris1999,naimanstone1998}. Degree theory can be used to obtain characterizations for structured classes of matrices~\cite[Chap. 6]{cottle1992lcp}. However, as of today, we lack a characterization of the degrees that can be realized by $Q$-matrices. We provide a characterization of this degree---and so a new sufficient criterion to be a $Q$-matrix---by studying the following problem, which is equivalent to the $\operatorname{LCP}(q,M)$ \cite{mangasarian2006absval}.
	Let $b\in\mathbb R^n$ and  $A\in\MnR$. Then the \emph{absolute value equation} (AVE) poses the problem to find a vector $z\in\mathbb R^n$ so that 
    \begin{equation}\label{eq:ave}
        z-A|z|\ =\ b\,,
    \end{equation}
    where $|\cdot|$ denotes the componentwise absolute value.
	The AVE is an interesting problem in its own right. For example, a result by Rump \cite[Thm. 2.8]{rump1997theorems} relates the number of solutions of \eqref{eq:ave} to the condition number of the matrix $A$, which is noteworthy in light of recent developments in real algebraic geometry that deal with precisely such connections of complexity and condition \cite[Part III]{buecuc:condition}.
	However, the main focus of theoretical investigations of the AVE is to obtain statements about the LCP. 
	A recent success of this approach is the development of condition numbers for the AVE that lead to new error bounds for the LCP~\cite{zamani2022condition}.
	
	We will study solvability of the AVE, but with our eyes on the $Q$-matrix problem. 
	To this end, following Cottle, Pang and Stone~\cite{cottle1992lcp}, we investigate the piecewise linear function
\begin{equation}\label{eq:affinemapFA}
\begin{aligned}
    F_A:\bbRn&\rightarrow \bbRn\\
    z&\mapsto z-A|z|
\end{aligned}
\end{equation}
	associated to the AVE \eqref{eq:ave} and determine its degree and its degree modulo 2 (see Section~\ref{def:degreemap}) in terms of the \emph{aligned spectrum} of $A$ (see Section~\ref{sec:alignedspectrum}):
	\begin{equation}\label{eq:formulaalignedspec}
	    \specA(A):=\{\lambda\geq 0\mid \exists x\neq 0\,:\,A|x|=\lambda x\},
	\end{equation}
	whose elements $\lambda\in \specA(A)$ are called the \emph{aligned values} of $A$.
	The first main result of this article relates the degree of $F_A$ to what we call the \emph{aligned count} of $A$:
	\begin{equation}\label{def:alignedcount}
	   \countA(A):=\#\{\lambda\in\specA(A)\mid \lambda>1\}
	\end{equation}
	where the count on the right-hand side is with multiplicities.
	
	We state now the main results of this paper, whose proofs are left to Section~\ref{sec:proof}. The term \emph{generic} in the theorem and the corollaries below means that the statement holds for all matrices that have a specific property (see Definition~\ref{defi:generic}), which is satisfied for all matrices except those in a given homogeneous hypersurface (see Section~\ref{sec:generic}). This condition, akin to the general position condition in the polyhedral world, guarantees that a random matrix (with respect to a continuous distribution) is generic with probability $1$.
	
    \begin{theorem}\label{thm:degswitch}
		Let $A\in\MnR$ be generic such that $1\notin \specA(A)$. Then the degree of $F_A$ is well-defined and it satisfies that
		\begin{equation}
		\deg F_A\ \equiv\ 1+\countA(A)\mod 2\,
		\end{equation}
		Moreover, $\deg F_A$ equals $1$ if all aligned values are smaller than $1$, and it equals $0$ if all aligned values are larger than $1$.  
	\end{theorem}
	\begin{corollary}\label{cor:numberalignedvalues}
	    Let $A\in\MnR$ be a generic matrix. Then the number of aligned values of $A$, counted with multiplicity, is odd.
	\end{corollary}
	\begin{corollary}\label{cor:solvabilityAVE}
	    Let $A\in\MnR$ be a generic matrix such that $1\not\in\specA(A)$. If $\countA(A)$ is even, then the AVE \eqref{eq:ave} has a solution for every $b\in\bbRn$.
	\end{corollary}


After stating Theorem~\ref{thm:degswitch}, we might wonder if there is an exact formula for the degree of $F_A$ when $A$ is generic (in the sense of Definition~\ref{defi:generic}). Indeed, there is such a formula and Theorem~\ref{thm:degswitch} is a direct consequence of the following more general---but more technical---theorem. 

\begin{theorem}\label{thm:exactformula}
Let $A\in\MnR$ be generic and such that $1\notin \specA(A)$. Then the degree of $F_A$ is well-defined and it satisfies that
\[
\deg F_A=1-\sum\{\sign(\chi_{SA}'(\lambda))\mid \lambda>1, S\in\mathcal{S}, \exists x\in \bbRn_{> 0}\,:\,SAx=\lambda x\}
\]
where $\chi_{SA}$ is the characteristic polynomial of $SA$, and $\mathcal{S}\subseteq\MnR$ is the set of sign matrices, i.e., diagonal matrices with $\pm1$ in the diagonal entries.
\end{theorem}

Observe that the right-hand side sum runs over all aligned values greater than 1, since we have that $A|x|=\lambda x$ for some $x\neq 0$ if and only if $SAx=x$ for some $S\in\Signs$ and $x\in\mathbb{R}^n_{\geq 0}$. Now, for a generic matrix (see Definition~\ref{defi:generic}), all aligned values correspond to simple eigenvalues of some $SA$, and so the right-hand side sum is nothing more than a ``signed aligned count'', i.e., a signed variation of the aligned count $\countA(A)$. In this way, Theorem~\ref{thm:degswitch} is just Theorem~\ref{thm:exactformula} reduced modulo 2.

As we mentioned above, a key reason to study the AVE is to gain insights into the equivalent LCP. 
To this end we derive LCP-analogues for the concept of the aligned spectrum and all statements about the degree of $F_A$ listed in this introduction.
Concerning our afore-stated interest in $Q$-matrices, we note that Corollary \ref{cor:solvabilityAVE} directly translates into a statement about the latter, i.e., the coefficient matrix of an LCP is a $Q$-matrix if the LCP-equivalent of the aligned count is even (Corollary \ref{cor:solvabilityLCP}). We observe, however, that this is only a sufficient condition for being a $Q$-matrix as shown by the example in \cite[p. 187--188]{kellywatson1979}, and that checking computationally this condition will probably be hard, just as in the case of the computation of the sign-spectral radius~\cite[Cor. 2.9]{rump1997theorems}.

\begin{remark}\label{remark:GAVE}
We note that our techniques extend straightforwardly to the case of \emph{generalized absolute value equations} (GAVE) given by
\[Bz-A|z|=b\]
with $A,B\in\MnR$ and $b\in\bbRn$. The reason for this is that our results are for generic matrices $A$ and $B$, and so we can assume, without loss of generality, that $B$ is invertible (a generic assumption) and reduce the generic GAVE above to the following generic AVE:
\[z-B^{-1}A|z|=b.\]
This same remark applies to the more general horizontal linearity complementarity problems which we discuss in Remark~\ref{remark:HLCP} at the end of Section~\ref{sec:transfertoLCPs}.
\end{remark}

\paragraph*{Organization}

In Section~\ref{sec:alignedspectrum}, we recall the sign-real spectrum and we introduce the aligned spectrum which naturally emerges during the study of the eigenproblem of $F_A$; in Section~\ref{sec:degree}, we recall the topological notion of degree, its formula in terms of signed counts of preimages of a regular value, and some of its properties in the case of interest; in Section~\ref{sec:generic}, we introduce the notion of genericity of a matrix relevant to our context and show some perturbation results. In Section~\ref{sec:transfertoLCPs}, we show how the above results for AVEs can be transferred to LCPs. We conclude in Section~\ref{sec:proof} proving Theorem~\ref{thm:degswitch} and its corollaries, and also Theorem~\ref{thm:exactformula}.

\paragraph*{Acknowledgements}

M.R. would like to thank Jiri Rohn and Siegfried M. Rump for helpful discussions.
Further, he would like to thank his advisor, Michael Joswig, for general support.

J.T.-C. was supported by a postdoctoral fellowship of the 2020 ``Interaction'' program of the \emph{Fondation Sciences Mathématiques de Paris}, and was partially supported by ANR JCJC
GALOP (ANR-17-CE40-0009), the PGMO grant ALMA, and the PHC GRAPE. He is also grateful to Evgenia Lagoda for moral support and Gato Suchen for the mathematical discussions regarding the proof of Theorem~\ref{thm:exactformula}.

We thank the referees for useful suggestions that helped to improve this paper.

\section{Sign-Real and Aligned Spectra}\label{sec:alignedspectrum}

The sign-real and aligned spectra of $A$ emerge naturally when we study the map $F_A$~\eqref{eq:affinemapFA}. Note that studying this map in order to understand \eqref{eq:ave} is similar to the strategy in linear algebra to study $z\mapsto Az$ in order to understand the solvability of $Az=b$

We note that $F_A$ is a positively homogeneous map, i.e., for every $z\in\bbRn$ and $\lambda>0$, $F_A(\lambda z)=\lambda F_A(z)$; and that $F_A$ is piecewise linear, with linear parts of the form
\begin{equation}
    \mathbb{I}-AS
\end{equation}
for sign matrices $S\in\Signs:=\{\mathrm{diag}(s_1,\ldots,s_n)\mid s_i\in\{-1,+1\}\}$.

\subsection{Sign-real spectrum and bijectivity}

The sign-real spectrum was used independently by Rump \cite{rump1997theorems} and---expressed in the language of interval arithmetic---by Rohn \cite[Thm. 5.1]{rohn1989interval} (cf.~\cite[Chap. 6]{neumaier1990interval}) to determine when the function $F_A$ is bijective, or equivalently, when the absolute value equation \eqref{eq:ave} is uniquely solvable for an arbitrary vector $b$. 

\begin{definition}
Let $A\in \MnR$. The \emph{sign-real spectrum} of $A$, $\specSR(A)$, is the set
\[
\specSR(A):=\bbR_{\geq 0}\cap\bigcup_{S\in\Signs}\spec(SA).
\]
\end{definition}

\begin{theorem}\label{theorem:solvabilityAVE}
Let $A\in\MnR$. Then the following are equivalent: 
\begin{enumerate}[(a)]
    \item $\max\specSR(A)<1$,
    \item $F_A$ is bijective,
    \item The AVE \eqref{eq:ave} has a unique solution for every $b\in\bbRn$.\eproof
\end{enumerate}
\end{theorem}

The largest element of the sign-real spectrum, $\max\specSR(A)$, is called the \emph{sign-real spectral radius}. Due to Theorem \ref{theorem:solvabilityAVE} it can be considered as a generalization of contractivity conditions for linear operators. Rump~\cite{rump1997theorems} also showed that the sign-real spectral radius generalizes the Perron root to matrices that are not positive. Another generalization of Perron Frobenius theory, to homogeneous monotone functions on the positive cone, was introduced in \cite{gaubert2004perron}. It is not clear how these two generalized theories are related apart from their common origin.

\subsection{Aligned spectrum and the eigenproblem}

The aligned spectrum arises when studying the eigenproblem for $F_A$: Determine for which $\lambda\geq 0$ and $v\in\bbRn\setminus 0$, we have
\[
F_A(v)=\lambda v.
\]

\begin{proposition}\label{prop:fixedpointsFA}
Let $A\in \MnR$, $\lambda\geq 0$ and $v\in\bbRn\setminus 0$. Then $F_A(v)=\lambda v$ if and only if there is some sign matrix $S\in\Signs$ such that $Sv\geq 0$ and $(1-\lambda,|v|)$ is an eigenpair of $SA$.
\end{proposition}
\begin{proof}
If $F_A(v)=\lambda v$, then we have that $(1-\lambda)v=A|v|$. Now, let $S\in \Signs$ such that $Sv\geq 0$, by taking as the diagonal elements of $S$ the signs of the components of $v$. Then $(1-\lambda)|v|=(1-\lambda)Sv=SA|v|$. Hence there is $S\in\Signs$ such that $Sv\geq 0$ and $(1-\lambda,|v|)$ is an eigenpair of $SA$.

Conversely, if such an $S$ exists, then
\[
F_A(v)=v-A|v|=v-S(SA)|v|=v-S(1-\lambda)|v|=\lambda v,
\]
where the second equality uses $S^2=\mathbb{I}$, the third one that $(1-\lambda,|v|)$ is an eigenpair of $SA$ (we are assuming this) and the fourth one that $v=S|v|$ (which follows from the assumed $Sv\geq 0$).  
\end{proof}

In view of the above proposition, the aligned spectrum is introduced. We note that the definition below is equivalent to that given in~\eqref{eq:formulaalignedspec}.
\begin{definition}\label{defi:alignedtrios}
Let $A\in \MnR$. An \emph{aligned trio} of $A$ is a triplet $(\lambda,S,v)\in \bbR_{\geq 0}\times\Signs\times (\SnMi\cap\bbRn_{\geq 0})$ such that
\[
SAv=\lambda v.
\]
Given such a trio, we call $(S,v)$ an \emph{aligned vector} and $\lambda$ an \emph{aligned value} of $A$. The \emph{aligned spectrum} of $A$, which we denote $\specA(A)$, is the set of aligned values of $A$, i.e.,
\[
\specA(A):=\{\lambda\geq 0\mid \exists S\in\Signs,\,v\in\bbRn_{\geq0}\setminus 0\,:\,SAv=\lambda v\}.
\]
\end{definition}

The following proposition shows how the aligned spectrum is related to the solution set of $F_A(z)=0$. In analogy to linear maps, we call $F_A$ \emph{nondegenerate} if $F_A(z)=0$ has only the trivial solution $z=0$. 

\begin{proposition}\label{prop:nondegenerateFA}
Let $A\in\MnR$. Then $F_A(z)=0$ has non-trivial solutions if and only if $1\in\specA(A)$.
\end{proposition}
\begin{proof}
By Proposition~\ref{prop:fixedpointsFA}, non-trivial solutions of $F_A(z)=0$ correspond to aligned trios of the form $(1,S,v)$. Hence, the claim follows.
\end{proof}

From the definitions it follows that
\begin{equation}
    \specA(A)\subseteq \specSR(A).
\end{equation}
However, this is not an equality in general.

\begin{example}
Let
\begin{equation}
	A\define \frac{2}{\sqrt{33}+3}\begin{pmatrix}
	1 & 2 \\ -1 &4
	\end{pmatrix}.
\end{equation}
One readily checks that
\[
\specA(A)=\left\{\frac{\sqrt{33}-3}{\sqrt{33}+3},\frac{4}{\sqrt{33}+3},\frac{6}{\sqrt{33}+3}\right\}\subsetneq \specSR(A)=\left\{\frac{\sqrt{33}-3}{\sqrt{33}+3},\frac{4}{\sqrt{33}+3},\frac{6}{\sqrt{33}+3},1\right\}.
\]
Moreover, by Theorem~\ref{theorem:solvabilityAVE} and Proposition~\ref{prop:nondegenerateFA}, this example shows that $F_A$ might not be bijective, even though it is nondegenerate. Furthermore, it shows that the largest aligned value and the sign-real spectral radius do not necessarily coincide. In light of Theorem \ref{thm:exactformula}, this demonstrates that we may have $\deg F_A=1$ without bijectivity.
\end{example}
We finish with the following example which shows that $A\mapsto \max\specA(A)$ is not continuous unlike $A\mapsto \max\specSR(A)$ which is~\cite[Corollary~2.5]{rump1997theorems}. However, in the generic case (see Definition~\ref{defi:generic}), we can recover continuity, since simple real eigenvalues cannot become complex and strictly positive vectors cannot become nonpositive under an arbitrarily small perturbation.

\begin{example}
Let $t$ lie in a sufficiently small neighborhood of $0$ and consider the following family of matrices:
\begin{equation}\label{eq:unstable}
	A_t\define \begin{pmatrix}
	1 & -0.5-t \\ 0.5 &0
	\end{pmatrix}\,.
\end{equation}
A straightforward calculation shows that $\specSR(A_t)$ is equal to
\[
\left\{\frac{1+\sqrt{-2t}}{2},\frac{1-\sqrt{-2t}}{2},\frac{\sqrt{2}\sqrt{1+t}-1}{2},\frac{1+\sqrt{2}\sqrt{1+t}}{2}\right\},
\]
if $t\leq 0$; and
\[
\left\{\frac{\sqrt{2}\sqrt{1+t}-1}{2},\frac{1+\sqrt{2}\sqrt{1+t}}{2}\right\}
\]
if $t>0$. Similarly, we can see that $\specA(A_t)$ is equal to
\[
\left\{\frac{1+\sqrt{-2t}}{2},\frac{1-\sqrt{-2t}}{2},\frac{\sqrt{2}\sqrt{1+t}-1}{2}\right\},
\]
if $t\leq 0$; and
\[
\left\{\frac{\sqrt{2}\sqrt{1+t}-1}{2}\right\}
\]
if $t>0$.

Hence, we have that
\[
\max\specSR(A_t)\, =\, \frac{1+\sqrt{2}\sqrt{1+t}}{2}\, \geq\, \max\specA(A_t)\, =\, \begin{cases}
\frac{1+\sqrt{-2t}}{2}\;,&\text{if }t\leq 0\;,\\
\frac{\sqrt{2}\sqrt{1+t}-1}{2}\;,&\text{if }t>0\;.
\end{cases}
\]
This shows that the maximum of the aligned spectrum is not continuous.
\end{example}

\section{Degree of a map}\label{sec:degree}

The degree of a continuous map $G:\SnMi\rightarrow\SnMi$ is a fundamental topological invariant that is preserved under homotopy. Intutively, we only have to think of the degree as the number of times that the map wraps $\SnMi$ around itself. Its formal definition is as follows.

\begin{definition}\cite[p.~134]{hatcher:topology}\label{def:degreemap}
Let $G:\SnMi\rightarrow\SnMi$ be a continuous map. The \emph{degree} of $G$, $\deg G$, is the unique integer $d$ such that the induced map $H_{n-1}(f):H_{n-1}(\SnMi)\rightarrow H_{n-1}(\SnMi)$ of homology groups is given by
\[x\mapsto dx\]
under the choice of any fixed isomorphism $H_{n-1}(\SnMi)\simeq \mathbb{Z}$.
\end{definition}

 Among the main properties of the degree, we have the following~\cite[p.~134]{hatcher:topology} (cf.~\cite[p. 98 ff]{ruiz2009deg}).

\begin{proposition}\label{prop:degreebasics}Let $G_0,G_1:\SnMi\rightarrow\SnMi$ be continuous maps. Then the following holds:
\begin{enumerate}[(1)]
    \item $\deg \mathrm{id}_{\SnMi}=1$ and $\deg (-\mathrm{id}_{\SnMi})=(-1)^n$.
    \item $\deg G_1\circ G_0=\deg G_1\deg G_0$.
    \item If there is a \emph{homotopy} between $G_0$ and $G_1$, that is, a continuous map $H:[0,1]\times \SnMi\rightarrow \SnMi$ such that for all $x\in\SnMi$, $H(0,x)=G_0(x)$ and $H(1,x)=G_1(x)$, then
    \[
    \deg G_0=\deg G_1.
    \]
    Moreover, the converse statement is also true.
    \item If $G_0$ is not surjective, then $\deg G_0=0$.\eproof
\end{enumerate}
\end{proposition}

Our investigation will be centered around $F_A$ which are nondegenerate. In this case, we can consider the spherical map
\begin{equation}
    \begin{aligned}
    \bF_A:\SnMi&\rightarrow \SnMi\\
    x&\mapsto F_A(x)/\|F_A(x)\|_2
\end{aligned}
\end{equation}
and define the \emph{degree of $F_A$} as
\begin{equation}
    \deg F_A:=\deg \bF_A.
\end{equation}
This definition agrees with a more traditional count used for maps $\bbRn\mapsto \bbRn$.
Recall that the \emph{set of regular values of $F_A$} is the set given by
\[
\mathrm{Reg}F_A:=\{y\in \bbRn\mid \forall x\in F_A^{-1}(y),\,\partial F_A(x)\text{ is well-defined and invertible}\},
\]
where $\partial F_A$ is the Jacobian of $F_A$.
\begin{proposition}
Let $A\in\MnR$ be such that $1\notin \specA(A)$. Then the set of regular values of $F_A$ is dense
and for all $y\in \mathrm{Reg}F_A$ we have
\begin{equation}
    \deg F_A=\sum_{x\in F_A^{-1}(y)}\sign(\det (\partial F_A(x))).
\end{equation}
\end{proposition}
\begin{proof}
By \cite[p. 509 ff]{cottle1992lcp} the oriented preimage counts of a nondegenerate positively homogeneous function and its restriction to the sphere coincide. 
\end{proof}

Moreover, the following proposition is helpful for our degree computations.

\begin{proposition}\label{prop:degreeFAbasics}
Let $A_0,A_1\in\MnR$. Then:
\begin{enumerate}[(1)]
    \item If $A:[0,1]\rightarrow\MnR$ is a continuous path between $A_0$ and $A_1$ such that for all $t\in [0,1]$, $1\notin\specA(A(t))$; then
    \[\deg F_{A_0}=\deg F_{A_1}.\]
    \item If $\specA(A)\subseteq [0,1)$, then $\deg F_{A}=1$. 
    \item If $\specA(A)\subset (1,\infty)$, then $\deg F_A=0$.
\end{enumerate}
\end{proposition}
\begin{proof}
(1) Consider the following homotopy between $\bF_{A_0}$ and $\bF_{A_1}$:
\[
[0,1]\times \SnMi\ni(t,x)\mapsto H(t,x)=\bF_{A(t)}(x)=\frac{F_{A(t)}(x)}{\|F_{A(t)}(x)\|_2}.
\]
This homotopy is well-defined because, by assumption and Proposition~\ref{prop:nondegenerateFA}, $F_{A(t)}(x)$ does not vanish at any $(t,x)$. Hence $\bF_{A_0}$ and $\bF_{A_1}$ are homotopic and so they have the same degree.

(2) Consider the path
\[
[0,1]\ni t\mapsto A(t):=(1-t)A.
\]
This path joins $A$ with $\mathbb{O}$ and it satisfies the condition of (1). Hence $\deg F_A=\deg F_{\mathbb{O}}=\deg \mathrm{id}_{\SnMi}=1$.

(3) Consider the following homotopy
\[
[0,1]\times \SnMi\ni(t,x)\mapsto H(t,x)=\frac{(1-t)x-A|x|}{\|(1-t)x-A|x|\|_2}.
\]
Note that the map $(t,x)\mapsto (1-t)x-A|x|$ is continuous, so if it is not vanishing, then the above homotopy is well-defined and continuous. If $t<1$, then
\[\SnMi\ni x\mapsto (1-t)x-A|x|=(1-t)F_{A/(1-t)}(x)\]
cannot vanish by Proposition~\ref{prop:nondegenerateFA}, since $1\notin \specA(A/(1-t))=\specA(A)(1-t)\subset (1,\infty)$. If $t=1$, then 
\[
\SnMi\ni x\mapsto A|x|
\]
does not vanish, because otherwise $0\notin \specA(A)$. Thus the desired map does not vanish and we obtain a homotopy between $\bF_A$ and
\[\SnMi\ni x\mapsto \frac{A|x|}{\|A|x|\|_2}.\]
If we precompose this map with $x\mapsto |x|$, it does not change. Now, since $x\mapsto |x|$ is not surjective,
\[
\deg \left(x\mapsto \frac{A|x|}{\|A|x|\|_2}\right)=\deg\left(x\mapsto \frac{A|x|}{\|A|x|\|_2}\right)\deg(x\mapsto |x|)=0,
\]
as we wanted to show.
\end{proof}

We conclude with an example which shows that the relationship between the aligned spectrum and the degree in Theorem~\ref{thm:degswitch} holds only modulo 2

\begin{example}\label{ex:modulotwounavoidable}
Let $\varepsilon$ be in a sufficiently small neigborhood of zero. Consider the family of matrices:
\begin{align}\label{eq:unstable2}
	B_\varepsilon\define \begin{pmatrix}
	2.5 & -1.25-\varepsilon \\ 1.25 &0
	\end{pmatrix}\,.
\end{align}
We can see that
\[
\specA(B_\varepsilon)=\left\{1.25\left(1+\sqrt{-0.8\varepsilon}\right),1.25\left(1-\sqrt{-0.8\varepsilon}\right),1.25\left(\sqrt{2}\sqrt{1+0.4\varepsilon}-1\right)\right\},
\]
if $\varepsilon \leq 0$; 
and that
\[
\specA(B_\varepsilon)=\left\{1.25\left(\sqrt{2}\sqrt{1+0.4\varepsilon}-1\right)\right\},
\]
if $\varepsilon>0$. And so we see that
\[
\#\{\lambda\in\specA(B_\varepsilon)\mid \lambda>1\}=\begin{cases}
2,&\text{if }\varepsilon<0\\
0,&\text{if }\varepsilon>0
\end{cases},
\]
and that
\[
\deg F_{B_\varepsilon}=1
\]
by Proposition~\ref{prop:degreeFAbasics}, cf. Figure \ref{plot:pert-ex}. 

On the one hand, this shows that the degree is more stable than the number of aligned values greater than one. On the other hand, note that the change in the number of aligned values greater than one happens because for $\varepsilon=0$, $B_\varepsilon$ is not generic---it has a double aligned value.

Moreover, for $\varepsilon<0$, $B_\varepsilon$ is generic (see Definition~\ref{defi:generic}) and it satisfies
\[
\#\{\lambda\in\specA(B_\varepsilon)\mid \lambda>1\}-1=\deg F_{B_\varepsilon};
\]
and for $\varepsilon>0$, $B_\varepsilon$ is still generic, but 
\[
\#\{\lambda\in\specA(B_\varepsilon)\mid \lambda>1\}+1=\deg F_{B_\varepsilon}.
\]
This shows that the equality modulo 2 in Theorem~\ref{thm:degswitch} cannot be corrected in an easy way to obtain an equality between the degree and the number of aligned values greater than one.
We need the more technical expression in Theorem \ref{thm:exactformula} for this. 

	\begin{figure}
		\includegraphics[width=0.47\textwidth]{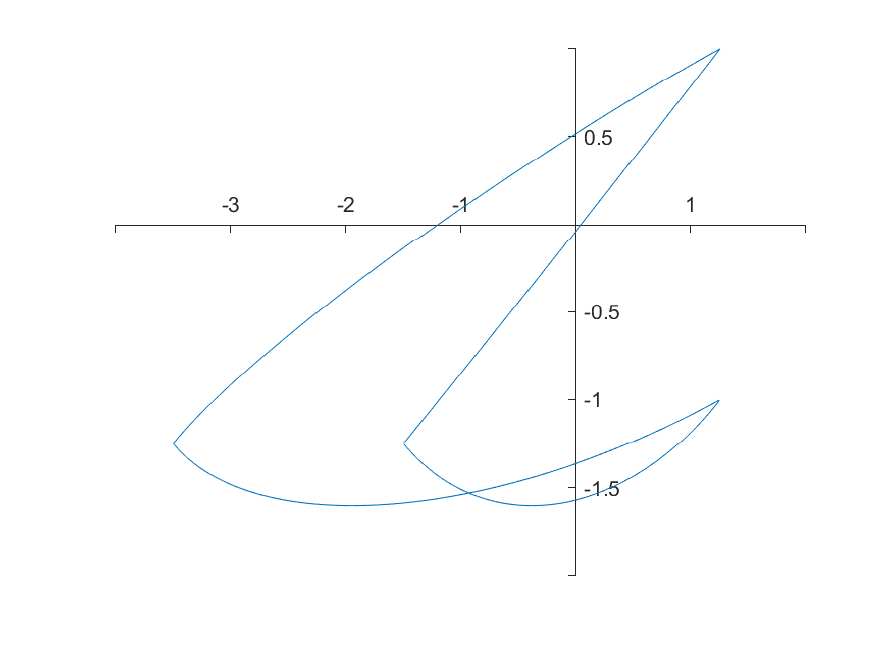}	
		\includegraphics[width=0.47\textwidth]{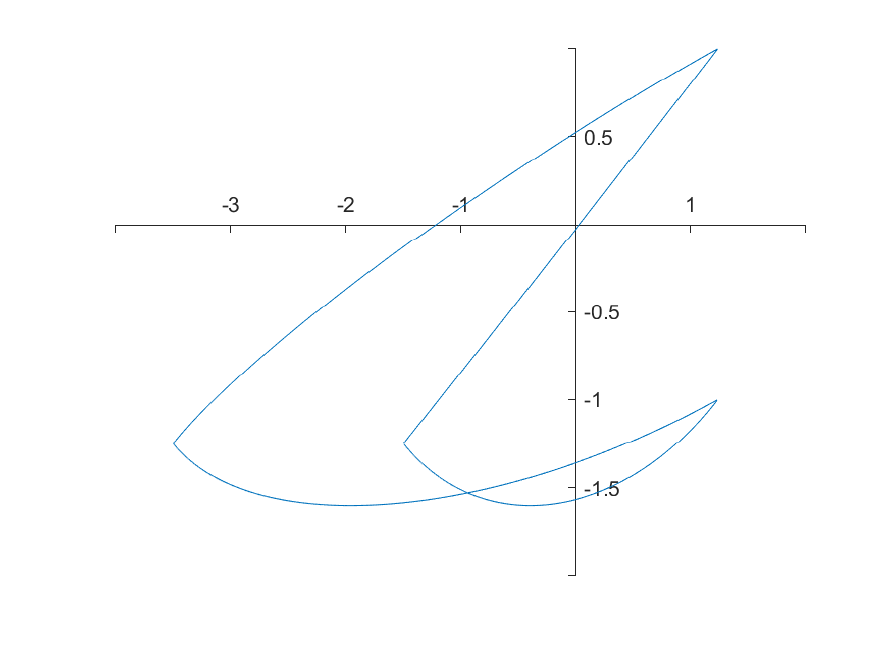}
		\caption{Image of the unit circle under $F_{B_\varepsilon}$. For $\varepsilon=-0.01$ (left) and $\varepsilon=0.01$ (right) it winds around the origin.}\label{plot:pert-ex}
	\end{figure}

\end{example} 

\section{Generic matrices}\label{sec:generic}

To make precise the statement of Theorem~\ref{thm:degswitch}, we introduce a notion of genericity adapted to our setting. 

\begin{definition}\label{defi:generic}
A \emph{generic} matrix is a matrix $A\in\MnR$ such that for every aligned trio $(\lambda,S,v)$ of $A$, a) $\lambda$ is a simple eigenvalue of $SA$, and b) $v$ is strictly positive.
\end{definition}

Since genericity usually means 
a class of entities whose complement is contained in a proper algebraic hypersurface, we need to show that the above notion is indeed a generic one according to the common use.

\begin{proposition}\label{prop:genericreallygeneric}
The set of matrices that are not generic in $\MnR$ is contained in a proper algebraic hypersurface.

In particular, for any random matrix $\mathfrak{A}\in\MnR$ with an absolutely continuous distribution, $\mathfrak{A}$ is generic almost surely.
\end{proposition}
\begin{proof}
We only have to prove that the set of matrices with double aligned values or aligned vectors in the boundary of the positive orthant is contained in a hypersurface.

The above set is contained in the union of the sets
\begin{equation}\label{eq:SetDiscA}
    \{A\in\MnR\mid S A\text{ has a double eigenvalue}\}
\end{equation}
and
\begin{equation}\label{eq:SetDiscB}
\{A\in\MnR\mid S A\text{ has a non-zero eigenvector in }H\}
\end{equation}
where $S\in\Signs$ runs over all sign matrices and $H$ over all coordinate hyperplanes---of the form $X_i=0$. Thus, if we show that each one of these sets is contained in a hypersurface, then we are done, since a finite union of hypersurfaces is a hypersurface.

On the one hand, 
the set~\eqref{eq:SetDiscA} is given by the discriminant of the characteristic polynomial of $SA$, which is well-known to define a proper algebraic hypersurface in $\MnR$.
On the other hand, the set~\eqref{eq:SetDiscB} is a proper algebraic hypersurface by \cite[Propoposition~1.2]{ottavianisturmfels2013}. Hence, all the sets are proper algebraic hypersurfaces, and the proof is complete.
\end{proof}

\subsection{Interpretation of count for generic matrices}
When we introduced the aligned count, see~\eqref{def:alignedcount}, 
\[\countA(A)\,=\,\#\{\lambda\in\specA(A)\mid \lambda>1\},\] 
we said that the right-hand side counts multiplicities. Note that for generic $A$, an aligned value $\lambda$ will always be a simple eigenvalue of the corresponding $SA$, where $S\in\Signs$. However, there might be more than one such $SA$. Because of this we have to include the \enquote{counted with multiplicity}. The following proposition gives an alternative interpretation of the central quantity $\countA(A)$ for Theorem~\ref{thm:degswitch} in terms of $\bF_A$.

\begin{proposition}\label{prop:genericcountig}
Let $A\in\MnR$ be generic. Then $\countA(A)$ is equal to the number of fixed points of $\bF_A$ that are images of their antipodal points under $\bF_A$. In other words,
\[
\countA(A)\,=\,\#\{x\in \SnMi\mid x=\bF_A(x)=\bF_A(-x)\}.
\]
\end{proposition}

For proving this proposition, the following proposition will be useful.

\begin{proposition}\label{prop:fixedpointsbFA}
Let $A\in\MnR$ and $x\in \SnMi$. If $F_A$ is nondegnerate, then $x$ is a fixed point of $\bF_A$ if and only if there is an aligned trio $(\lambda,S,v)$ such that either $x=-Sv$ or $\lambda<1$ and $x=Sv$. Moreover, when $x$ is a fixed point of $\bF_A$, the following are equivalent:
\begin{itemize}
    \item $\bF_A(-x)=x$.
    \item There is an aligned trio $(\lambda,S,v)$ such that $\lambda>1$ and $x=-Sv$.
\end{itemize}
\end{proposition}
\begin{proof}[Proof of Proposition~\ref{prop:genericcountig}]
    By Proposition~\ref{prop:fixedpointsbFA}, the fixed points $x\in\SnMi$ of $\bF_A$ such that $\bF_A(-x)=x$ are in one-to-one correspondence with the aligned trios $(\lambda,S,v)$ such that $\lambda>1$. Since $A$ is generic, this means precisely the number of aligned values greater than one counted with multiplicity.
\end{proof}
\begin{proof}[Proof of Proposition~\ref{prop:fixedpointsbFA}]
If $(\lambda,S,v)$ is an aligned trio, then
\[F_A(Sv)=(1-\lambda)Sv\]
and
\[
F_A(-Sv)=(1+\lambda)(-Sv).
\]
In this way, $-Sv$ is always a fixed point of $\bF_A$ and $Sv$ is so if and only if $\lambda<1$. This shows one direction.

If $x\in\SnMi$ is a fixed point of $\bF_A$, then for some $\mu>0$,
\[F_A(x)=\mu x.\]
Let $S\in \Signs$ be such that $Sx\geq 0$, so that $v=Sx$. Then we have that
\[
SAv=(1-\mu)v.
\]
If $1-\mu\geq 0$, then $(1-\mu,S,v)$ is an aligned trio such that $1-\mu<1$ and $x=Sv$. Otherwise, $1-\mu<0$, and then $(\mu-1,-S,v)$ is an aligned trio and $x=-(-S)v$. Hence there is an aligned trio $(\lambda,S,v)$ such that either $x=-Sv$ or $\lambda<1$ and $x=Sv$.

We show the second equivalence. Let $x\in\SnMi$ be a fixed point of $\bF_A$. Then, by the first part, there is an aligned trio $(\lambda,S,v)$ such that either $x=-Sv$ or $\lambda<1$ and $x=Sv$. In the second case, we have that
\[
\bF_A(-x)=-x.
\]
Thus we must have the first case. But then $\bF_A(-x)=x$ if and only if $\lambda>1$, because otherwise $-x$ is also a fixed point.

Suppose $x=-Sv$ for some aligned trio $(\lambda, S,v)$ such that $\lambda>1$. Then, by the first equivalence, $x=-Sv$ is a fixed point of $\bF_A$, and, by direct computation, $\bF_A(-x)=x$. 
\end{proof}

\subsection{Perturbation of matrices to make them generic}

The following proposition shows that matrices corresponding to nondegenerate maps can be slightly perturbed to obtain a generic matrix with the same corresponding degree.

\begin{proposition}\label{prop:randomperturbation}
Let $A\in\MnR$ be such that $1\notin \specA(A)$. Then:
\begin{enumerate}[(a)]
\item The quantity
\begin{equation}
\cond(A)\,\define\,\sup_{x\neq 0}\frac{\|x\|_2}{\|F_A(x)\|_2}
\end{equation}
is finite.
\item For every $\varepsilon \in \left(0,1/\cond(A)\right)$ and
\[
\tilde{A}\in B_F(A,\varepsilon):=\{X\in\MnR\mid \|X-A\|_F<\varepsilon\},
\]
$F_{\tilde{A}}$ is nondegenerate, and $\deg F_{\tilde{A}}=\deg F_A$.
\item Let $\tilde{A}\in B_F(A,\varepsilon)$ be a random matrix with the uniform distribution on $B_F(A,\varepsilon)$, then $\tilde{A}$ is generic with probability one.
\end{enumerate}
\end{proposition}
\begin{proof}
(a) We have that
\[
\min_{x\neq 0}\frac{\|F_A(x)\|_2}{\|x\|_2}
\]
is zero if and only if $1\in\specA(A)$ by Proposition~\ref{prop:nondegenerateFA}. Hence,
it is a positive number if $1\notin\specA(A)$ and its inverse, the quantity $\cond(A)$, must be finite.

(b) By the inequalities between matrix norms, we can show that
\[
\frac{1}{\cond(\tilde{A})}\,\geq\, \frac{1}{\cond(A)}-\left\|\tilde{A}-A\right\|_2\,\geq\,\frac{1}{\cond(A)}-\left\|\tilde{A}-A\right\|_F \,.
\]
Hence, if $\tilde{A}\in B_F(A,\varepsilon)$, with the given choice of $\varepsilon$, then no matrix in the segment $[A,\tilde{A}]$ can have $1$ as an aligned value. Consequently, we have a path between $A$ and $\tilde{A}$, given by $t\mapsto (1-t)A+t\tilde{A}$, such that no matrix in the path has $1$ as an aligned value, and so $\deg F_{\tilde{A}}=\deg {F_A}$ by Proposition~\ref{prop:degreeFAbasics} (1).

(c) This is a direct consequence of Proposition~\ref{prop:genericreallygeneric}.
\end{proof}

We observe that the above proposition can only be applied when $1\notin\specA(A)$. In that case, it allows us to produce a generic matrix $\tilde{A}$ such that $F_{\tilde{A}}$ has the same topological structure---degree---as $F_A$. Now, if $1\in\specA(A)$, this perturbation trick will not produce $F_{\tilde{A}}$ with the same topological structure as the following examples show.

\begin{example}
Let
\begin{equation}
A\define \begin{pmatrix}
		2 & -1 \\ -1 &0
		\end{pmatrix}
\end{equation}
for which $\specA(A)=\{1,\sqrt{2}-1\}$. Now consider the following perturbation
\begin{equation}
A_\varepsilon\define \begin{pmatrix}
		2 & -1-\varepsilon \\ -1 &0
		\end{pmatrix}.
\end{equation}
We can see that for $\varepsilon>0$, $\deg F_{A_\varepsilon}=1$, since all aligned values are smaller than one; and that for $\varepsilon<0$, $\deg F_{A_\varepsilon}=0$, after a straightforward computation. Moreover, Figure \ref{plot:circles} indicates---and it can be checked by computation---that
\[
F_{A_\varepsilon}(x)=\begin{pmatrix}r\\0
\end{pmatrix}\,,
\]
where $r>0$, does not have a solution for $\varepsilon<0$, cf. Figure \ref{plot:circles}.

Hence, when we perturb $A\in\MnR$ with $1\in\specA(A)$, neither do we get consistent topological information about $F_{\tilde{A}}$ via the perturbed matrix $\tilde{A}$, nor do we obtain consistent information about the general solvability of the AVE~\eqref{eq:ave}.

	\begin{figure}
		\includegraphics[width=0.47\textwidth]{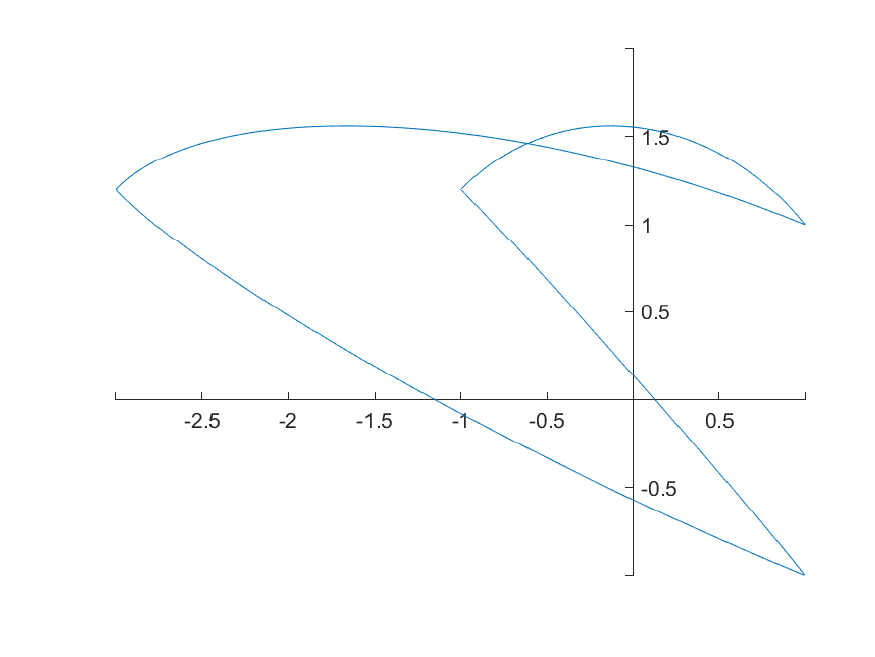}	
		\includegraphics[width=0.47\textwidth]{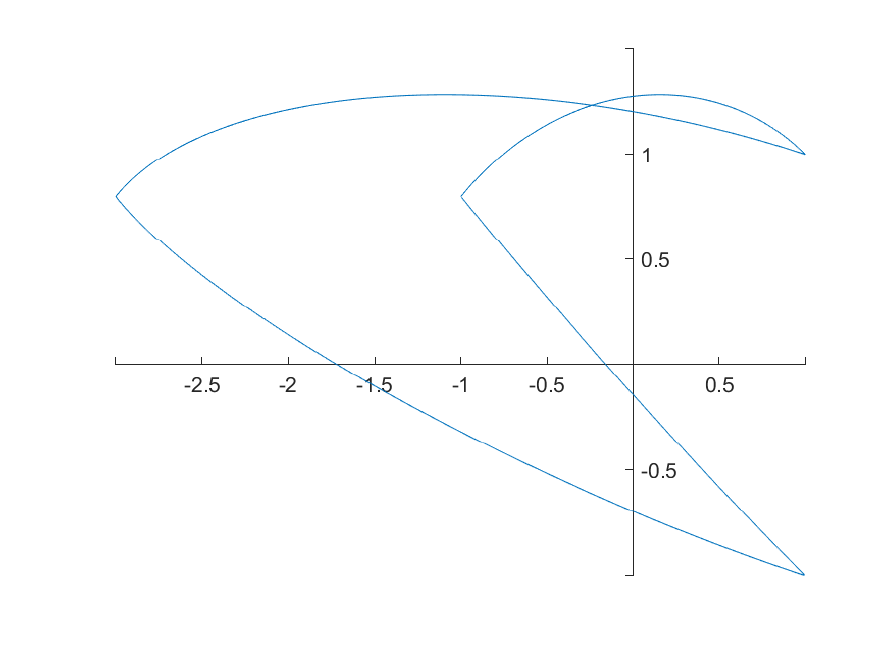}
		\caption{Image of the unit circle under $F_{A_\varepsilon}$. For $\varepsilon=0.2$ it winds around the origin (left). For $\varepsilon=-0.2$ it does not (right).}\label{plot:circles}
	\end{figure}

\end{example}	

\begin{example}
Consider a generic matrix $A\in\MnR$ so that all aligned values have odd multiplicity. Then for every $t>0$, $tA$ is generic as long as $1/t\notin \specA(A)$. As $t$ increases from zero to infinity, we have that $\deg F_{tA}$ alternates parity as $t$ crosses the inverses of the aligned values of $A$, by Theorem \ref{thm:degswitch}. This shows again that perturbing matrices with 1 as an aligned value will not produce consistent topological information.
\end{example}

\section{Transfer of results to LCPs}\label{sec:transfertoLCPs}

Classically~\cite[Chap. 1]{cottle1992lcp}, the LCP is associated to the following piecewise linear function  
\begin{equation}
	\begin{aligned}
		G_M:\bbRn&\rightarrow \bbRn\\
		z&\mapsto (M+\mathbb{I})z/2-(M-\mathbb{I})|z|/2\,, 
	\end{aligned}
\end{equation}
where we have rewritten the expression in terms of absolute values. If $G_M$ is surjective, then $\operatorname{LCP}(q,M)$ has a solution for any $q\in \bbRn$. Therefore studying the degree of $G_M$ is the $LCP$-equivalent of studying the degree of $F_A$ in the context of AVEs. 

The following proposition which is well-known in the literature (see~\cite[Prop. 1.4.4]{cottle1992lcp}, \cite{mangasarian2006absval} and \cite[Chap. 6]{neumaier1990interval}) makes explicit how AVEs, LCPs and $G_M$ relate.

\begin{proposition}\label{prop:LCPtoAVE}
	Let $M\in \MnR$ and $q\in\bbRn$. The map $(v,w)\mapsto v-w$ is a bijective correspondence between the solutions $(v,w)\in\bbRn_{\geq 0}\times\bbRn_{\geq 0}$ of $\operatorname{LCP}(q,M)$ and the solutions $z\in\bbRn$ of $G_M(z)=q$. In particular, if $M+\mathbb{I}$ is invertible, then $\operatorname{LCP}(q,M)$ is equivalent to solve the AVE given by
	\begin{equation*}\tag*{\qed}
		z-(M+\mathbb{I})^{-1}(M-\mathbb{I})|z|=2(M+\mathbb{I})^{-1}q
	\end{equation*}
\end{proposition}

In this case, the images of the spherical restrictions $\bF_{(M+\mathbb{I})^{-1}(M-\mathbb{I})}$ and $\bar{G}_M=G_M/\lVert G_M\rVert_2$ are congruent, which shows that the degrees of both maps are identical up to a sign change according to the determinant sign of the transformation matrix $M+\mathbb{I}$. In particular, the parities of their degrees are the same. 

To further analyze the degree of $G_M$, we introduce the LCP-variant of aligned trios of Definition~\ref{defi:alignedtrios} and of generic matrices of Definition~\ref{defi:generic}---for which we can prove claims analogous to those of Section~\ref{sec:generic}.
\begin{definition}
	Let $M\in \MnR$. An \emph{LCP-aligned trio} of $M$ is a triplet $(\lambda,S,v)\in \bbR_{\geq 0}\times\Signs\times (\SnMi\cap\bbRn_{\geq 0})$ such that
	\[
	(M-\mathbb{I})v\, =\, \lambda(M+\mathbb{I})Sv\,.
	\]
	Given such a trio, we call $(S,v)$ an \emph{LCP-aligned vector} and $\lambda$ an \emph{LCP-aligned value} of $M$. The \emph{LCP-aligned spectrum} of $M$, which we denote $\specLCP(M)$, is the set of LCP-aligned values of $M$. 
\end{definition}
\begin{remark}
In many of our results, we assume that $1\not\in\specLCP(M)$. This is equivalent to assuming that $M$ is an $R_0$-matrix, that is, the $\LCP(0,M)$ has a unique solution~\cite[Def.~3.8.7]{cottle1992lcp}. The latter is a common assumption when working with LCPs.
\end{remark}

LCP-aligned vectors are not eigenvectors of $G_M$ and thus also not fixed points of $\bar{G}_M$. However, they and their polar opposites are exactly the pairs of antipodal points which are mapped to pairs of antipodal points by $\bar{G}_M$ (if the corresponding LCP-aligned value is smaller than $1$) or to a single point (if the corresponding LCP-aligned value is larger than $1$).
This property of LCP-aligned vectors seems to be more crucial to the parity of the degree, which mirrors the number of times that the sphere is folded onto itself by $\bar{G}_M$, 
than the property of being an eigenvector or a fixed point of the spherical map.

\begin{definition}
	An LCP-generic matrix $M\in\MnR$ is a matrix such that a) $M+\mathbb{I}$ is invertible, and b)
	$(M+\mathbb{I})^{-1}(M-\mathbb{I})$ is generic.
\end{definition}
LCP-generic matrices are really generic due to Proposition \ref{prop:genericreallygeneric} and the fact that $\{M\in\MnR\mid \det(M+\mathbb{I})=0\}$ is an algebraic hypersurface.
Let $M\in\MnR$ be a matrix so that $G_M$ is nondegenerate. Analogously to Proposition \ref{prop:randomperturbation}, a random perturbation of $M$ will be LCP-generic almost surely.

We can now state the main theorems (Theorems~\ref{theorem:solvabilityAVE} and \ref{thm:exactformula}) and their corollaries in the context of LCPs. Recall that $\deg G_M$ is the degree of the spherical map $\bar{G}_M:x\mapsto G_M(x)/\|G_M(x)\|_2$. 

\begin{theorem}\label{thm:LCPmain1}
	Let $M\in\MnR$ be LCP-generic such that $1\not\in\specLCP(M)$.
	Then the degree of $G_M$ is well-defined and it satisfies that
	\[
	\deg G_M\,\equiv\, 1+
	\countLCP(M)\mod 2\,
	\]
	where $\countLCP(M)$ is the number of LCP-aligned values greater than one, counted with multiplicity. Moreover, we have $\deg G_M=\sign(\det(M+\mathbb{I}))$ if $\countLCP(M)=0$ and $\deg G_M=0$ if all LCP-aligned values of $M$ are larger than $1$.
\end{theorem}

\begin{theorem}\label{thm:LCPmain2}
	Let $M\in\MnR$ be LCP-generic such that $1\not\in\specLCP(M)$ 
	Then the degree of $G_M$ is well-defined and it satisfies that
	\[
	\deg G_M\,=\, \sign(\det(M+\mathbb{I})) -\sum_{\substack{(\lambda,S,v)\text{ LCP-aligned trio of }M\\ \lambda>1}} \sign\left(\chi_{M+\mathbb{I},(M-\mathbb{I})S}'(\lambda)\right)\,.
	\]
	where $\chi_{U,V}:=\det(TU-V)$ is the generalized characteristic polynomial of $(U,V)\in\MnR^2$.
\end{theorem}

\begin{corollary}\label{cor:numberLCPalignedvalues}
	Let $M\in\MnR$ be an LCP-generic matrix. Then the number of LCP-aligned values of $M$, counted with multiplicity, is odd.\eproof
\end{corollary}
\begin{corollary}\label{cor:solvabilityLCP}
	Let $M\in\MnR$ be LCP-generic such that $1\not\in\specLCP(M)$. If $\countLCP(M)$ is even, then $M$ is a $Q$-matrix, i.e., for all $q\in\bbRn$, $\operatorname{LCP}(q,M)$ has a solution.
\end{corollary}

To show these theorems and their corollaries, we only need to show how the notions in the LCP-setting crrespond to the notions in the AVE-setting. The following three propositions allow us to do these translations. The first one gives how aligned trios correspond to LCP-aligned trios; the second one gives sufficient condition for the degree of $G_M$ to be well-defined; and the third one relates the degree of $G_M$ to that of $F_A$ for an appropriate $A$.

\begin{proposition}\label{prop:spectra-correspondence}
	Let $M\in\MnR$ be such that $M+\mathbb{I}$ is invertible.
	Then $(\lambda,S,v)\in \bbR_{\geq 0}\times\Signs\times (\SnMi\cap\bbRn_{\geq 0})$ is an LCP-aligned trio of $M$ if and only if it is an aligned trio of $(M+\mathbb{I})^{-1}(M-\mathbb{I})$.
\end{proposition}

\begin{proposition}\label{prop:nondegenerateGM}
	Let $M\in\MnR$. Then $G_M(z)=0$ has non-trivial solutions if and only if $1\not\in\specLCP(M)$. In particular, if $1\not\in\specLCP(M)$, the degree of $G_M$ is well-defined, the set of regular values of $G_M$ is dense
	and for all $y\in \mathrm{Reg}G_M$ we have
	\begin{equation}
		\deg G_M=\sum_{x\in G_M^{-1}(y)}\sign(\det (\partial G_M(x))).
	\end{equation}
\end{proposition}

\begin{proposition}\label{prop:degree-correspondence}
	Let $M\in\MnR$ be such that $1\not\in\specLCP(M)$ and such that $M+\mathbb{I}$ is invertible. Then 
	\[
	\deg G_M\,=\,\sign(\det(M+\mathbb{I}))\cdot\deg F_{(M+\mathbb{I})^{-1}(M-\mathbb{I})}.
	\]
\end{proposition}

\begin{proof}[Proof of Proposition~\ref{prop:spectra-correspondence}]
	Note that, when $M+\mathbb{I}$ is invertible, $(\lambda,S,v)$ is a LCP-aligned trio of $M$ if and only if $\lambda Sv=(M+\mathbb{I})^{-1}(M-\mathbb{I})v$. The latter is equivalent to $\lambda v=S(M+\mathbb{I})^{-1}(M-\mathbb{I})v$, which means exactly that $(\lambda,S,v)$ is an aligned trio of $(M+\mathbb{I})^{-1}(M-\mathbb{I})$.
\end{proof}
\begin{proof}[Proof of Proposition~\ref{prop:nondegenerateGM}]
	The first claim is proven as in Proposition~\ref{prop:nondegenerateFA}. The second one follow from \cite[p. 509 ff]{cottle1992lcp}, since $G_M$ is a nondegenerate positively homogeneous map.
\end{proof}
\begin{proof}[Proof of Proposition~\ref{prop:degree-correspondence}]
	This follows from the multiplicative property of the degree and the fact that
	\[L_M\circ\bF_{(M+\mathbb{I})^{-1}(M-\mathbb{I})}=\bar{G}_M\]
	where $L_M:x\mapsto (M+\mathbb{I})x/\|(M+\mathbb{I})x\|_2$. Recall that we have that $\deg L_M=\sign(\det(M+\mathbb{I}))$.
\end{proof}

We now give the proof of Theorems~\ref{thm:LCPmain1} and~\ref{thm:LCPmain2}.

\begin{proof}[Proof of Theorem~\ref{thm:LCPmain1}]
	By Propositions~\ref{prop:spectra-correspondence}, \ref{prop:nondegenerateGM} and \ref{prop:degree-correspondence}, we have that
	\[
	\deg G_M\,\equiv\, \deg F_{(M+\mathbb{I})^{-1}(M-\mathbb{I})}\mod 2\,
	\]
	and
	\[
	\countLCP(M)\,\equiv\,\countA((M+\mathbb{I})^{-1}(M-\mathbb{I}))\mod 2.
	\]
	Hence the theorem follows by Theorem~\ref{thm:degswitch}.
\end{proof}
\begin{proof}[Proof of Theorem~\ref{thm:LCPmain2}]
	By Propositions~\ref{prop:spectra-correspondence}, 
	\[
	\deg G_M\,=\, \sign(\det(M+\mathbb{I}))\deg F_{(M+\mathbb{I})^{-1}(M-\mathbb{I})}.
	\]
	Hence, by Theorem~\ref{thm:exactformula} and Proposition~\ref{prop:degree-correspondence},
	\[
	\deg G_M=\sign(\det(M+\mathbb{I}))-\sum_{\substack{(\lambda,S,v)\text{ LCP-aligned trio of }M\\ \lambda>1}} \sign\left(\det(M+\mathbb{I})\chi_{(M+\mathbb{I})^{-1}(M-\mathbb{I})S}'(\lambda)\right).
	\]
	Now,
	\[
	\det(M+\mathbb{I}))\chi_{(M+\mathbb{I})^{-1}(M-\mathbb{I})S}'(\lambda)=\chi_{M+\mathbb{I},(M-\mathbb{I})S}'(\lambda),
	\]
	so we obtain the desired result.
\end{proof}

\begin{remark}\label{remark:HLCP}
Following with Remark \ref{remark:GAVE}, we can apply the result above to the more general \emph{horizontal linear complementarity problem} (HLCP) which, given $M,N\in\MnR$ and $q\in \bbRn$, asks to find $v,w\in\bbRn_{\geq 0}$ with $v^Tw=0$ so that 
\[
Nv - Mw = q\,.
\]
In the same way we show in this section that the LCP and the AVE are equivalent, we can show that this HLCP is equivalent to the following GAVE
\[
\frac{1}{2}(M+N)z-\frac{1}{2}(M-N)|z|=q.
\]
Then, as our results hold for generic matrices, we can translate the results in the particular setting of LCPs to HLCPs. Moreover, for generic (and thus invertible) $N$, the HCLP above is equivalent to $\LCP(N^{-1}q,N^{-1}M)$.
\end{remark}

\section{Proof of Theorems~\ref{thm:degswitch} and~\ref{thm:exactformula} and Corollaries~\ref{cor:numberalignedvalues} and \ref{cor:solvabilityAVE}}\label{sec:proof}

We first show how Theorem~\ref{thm:degswitch} follows from Theorem~\ref{thm:exactformula}. Then we prove the corollaries~\ref{cor:numberalignedvalues} and \ref{cor:solvabilityAVE}. We finish giving the proof of Theorem~\ref{thm:exactformula}.

\subsection{Proof of Theorem~\ref{thm:degswitch}}

Note that the formula of Theorem~\ref{thm:degswitch} can be rewritten as
\[
\deg F_A=1-\sum_{\substack{(\lambda,S,v)\text{ aligned trio of }A\\\lambda>1}}\sign(\chi_{SA}'(\lambda)).
\]
Now, since $A$ is generic, we have that for each aligned trio $(\lambda,S,A)$, $\lambda$ is a simple eigenvalue of $SA$ and so a simple root of $\chi_{SA}$. Hence $\sign(\chi_{SA}'(\lambda))$ is either $+1$ or $-1$ for each summand in the sum. Moreover, for a specific aligned value $\lambda$, the summand $\sign(\chi_{SA}'(\lambda))$ appears as many times as the size of
\[
\{(S,v)\mid S\in\Signs,\,v\in\SnMi\cap\bbR^n_{\geq 0}, (S,v)\text{ is the aligned vector of }A\text{ corresponding to }\lambda\}.
\]
But this quantity, for generic $A$, is precisely the multiplicity of $\lambda$. Hence, we have that for all $\lambda>1$,
\[
\sum_{\substack{(\lambda,S,v)\text{ aligned trio of }A}}\sign(\chi_{SA}'(\lambda))
\]
is the multiplicity of $\lambda$ mod 2. Hence
\[
\sum_{\substack{(\lambda,S,v)\text{ aligned trio of }A\\\lambda>1}}\sign(\chi_{SA}'(\lambda))\equiv \countA(A)\mod 2,
\]
and the first part of the theorem follows.

The second part is just Proposition~\ref{prop:degreeFAbasics} (2)-(3).

\subsection{Proof of Corollary~\ref{cor:numberalignedvalues}}

Since $A$ is generic, the multiplicity of $0$ as an aligned value is given by the size of the set
\[
\{(S,v)\mid S\in\Signs,\,v\in\SnMi\cap\bbR^n_{\geq 0}, (S,v)\text{ is the aligned vector of }A\text{ corresponding to }0\}.
\]
This set is invariant under the transformation $(S,v)\mapsto (-S,v)$. Hence, the multiplicity of $0$ as an aligned value is always even. 

By the proof of Lemma~\ref{lem:auxlemma}, we can consider perturbed matrix $\tilde{A}$ with the same number of positive aligned values as $A$, but such that $0$ is not an aligned value of $\tilde{A}$. Since the multiplicity of $0$ is even, this means that the parity of the number of aligned values of $\tilde{A}$ and $A$ is the same. Therefore, without loss of generality, we can assume that $A$ has only positive aligned values.

For $t>0$, we have that
\begin{equation}\label{eq:dilationaspectrum}
    \specA(tA)=t\specA(A).
\end{equation}
For some $t>0$ sufficiently large, every aligned value of $tA$ is larger than one, because, by assumption, $A$ has only positive aligned values. Thus, by the second part of Theorem~\ref{thm:degswitch}, $\deg F_{tA}=0$, and so, by the first part of the same theorem, $\countA(tA)$ is odd. Hence $tA$ has an odd number of aligned values, as we wanted to show.

\subsection{Proof of Corollary~\ref{cor:solvabilityAVE}}

If $\countA(A)$ is even, then, by Theorem~\ref{thm:degswitch}, $\deg F_A$ is odd, and so, in particular, non-zero. Hence, by Proposition~\ref{prop:degreebasics}(4), $\bF_A$ is surjective, and so is $F_A$. 

\subsection{Proof of Theorem~\ref{thm:exactformula}}

By perturbing $A$ randomly, we can assume without loss of generality that all aligned values of $A$ are simple, i.e., that every aligned value of $A$ appears only in one aligned trio $(\lambda,S,v)$ of $A$. The following lemma allows us to do this.

\begin{lemma}\label{lem:auxlemma}
Let $A\in\MnR$ be generic such that $1\notin\specA(A)$. Then there is $\tilde{A}$ such that 1) $1\notin \specA(A)$, 2) $\tilde{A}$ is generic, 3) $F_A$ and $F_{\tilde{A}}$ have the same degree, 4) $A$ and $\tilde{A}$ have the same aligned count, i.e., $\countA(A)=\countA(\tilde{A})$; and 
\[
\sum_{\substack{(\tilde{\lambda},S,v)\text{ aligned trio of }\tilde{A}\\\lambda>1}}\sign(\chi_{S\tilde{A}}'(\tilde{\lambda}))=\sum_{\substack{(\lambda,S,v)\text{ aligned trio of }A\\\lambda>1}}\sign(\chi_{SA}'(\lambda))
;\]
and 5) for every aligned value $\tilde{\lambda}$ of $\tilde{A}$, there is a unique $S\in\mathcal{S}$ such that $\tilde{\lambda}$ is an eigenvalue of $S\tilde{A}$.
\end{lemma}
\begin{proof}
By Proposition~\ref{prop:randomperturbation}, 1), 2), and 3) are guaranteed by taking $\tilde{A}$ in a sufficiently small neighborhood $B_F(A,\varepsilon)$ of $A$. Now, since $A$ is generic, for each $S\in\Signs$, given an eigenvalue $\lambda$ of $SA$, we proceed as follows:
\begin{enumerate}[(a)]
    \item If $\lambda\notin [1,\infty)\subseteq \mathbb{C}$, then, by continuity of the eigenvalues, we have that for an arbitrarily small perturbation of $A$, $\tilde{A}$, the corresponding eigenvalue, $\tilde{\lambda}$, is still outside $[1,\infty)$.
    \item If $\lambda>1$ is not an aligned value, then $(\lambda\mathbb{I}-SA)v$ does not vanish for $v\in\SnMi\cap\bbR^n_{\geq 0}$. Therefore
    \[
    \min_{v\in\SnMi\cap\bbR^n_{\geq 0}}\|(\lambda\mathbb{I}-SA)v\|_2>0.
    \]
    But this quantity is not only continuous in $\lambda$ and $A$, but $1$-Lipschitz in them. Hence, for an arbitrarily small perturbation $\tilde{A}$ of $A$, we can guarantee that the corresponding eigenvalue $\tilde{\lambda}$ does not become an aligned value of $\tilde{A}$.
    \item If $\lambda>1$ is an aligned value, consider an aligned trio $(\lambda,S,v)$ such that $\|v\|_2=1$. Since $A$ is generic, $\lambda$ is simple, and so we can apply the implicit function theorem to
    \[
    \bbR\times\SnMi\times\MnR\ni (\lambda,x,M)\mapsto (\det(\lambda\mathbb{I}-MS),(\lambda\mathbb{I}-MS)v)
    \]
    at $(\lambda,v,A)$. Hence, in a small neighborhood of $A$, we can write $\lambda$ and $v$ as smooth functions of $A$. Since $v$ is strictly positive, this means that for an arbitrarily small perturbation $\tilde{A}$ the aligned value $\lambda$ goes to an aligned value $\tilde{\lambda}$ which is still simple as an eigenvalue of $\tilde{A}S$. Moreover, if the perturbation is sufficiently small, $\tilde{\lambda}$ remains inside $(1,\infty)$, by continuity, and the signs of $\chi_{SA}'(\lambda)$ and of $\chi'_{S\tilde{A}}(\tilde{\lambda})$ coincide, by the continuity of $\chi'_{SM}(\mu)$ with respect to $(M,\mu)$.
\end{enumerate}
Putting the above together, we have that taking $\tilde{A}$ in a sufficiently small neighborhood $B_F(A,\varepsilon)$ of $A$ we can guarantee that it satisfies 1), 2), 3) and 4).

For 5), we only have to show that for almost all $\tilde{A}\in\MnR$, the $S\tilde{A}$, with $S\in\Signs$, do not share any eigenvalue. Once this is done, we can guarantee, by the above, that we can choose a perturbation $\tilde{A}$ with the desired properties. Note that this is the only point where the perturbation is not arbitrary.

We will show that the set 
\[M_S:=\{
X\in\MnR\mid X\text{ and }SX\text{ have an eigenvalue in common}
\}\]
is a proper algebraic hypersurface. Then the set of $X$ such that $S_1X$ and $S_2X$ share an eigenvalue for some $S_1,S_2\in \Signs$ is given by
\[
\bigcup_{S,T\in\Signs}SM_T.
\]
Therefore it will be a proper algebraic hypersurface, since it is a finite union of proper algebraic hypersurfaces. Hence, we can choose $\tilde{A}\in B_F(A,\varepsilon)$ such that $\tilde{A}$ does not lie in it.

The determinant of the Sylvester matrix of the characteristic polynomials of $X$ and $SX$ is zero if and only if $X$ and $SX$ share an eigenvalue (see~\cite[Ch. 3, Prop. 8]{coxbook}). Hence $M_S$ is described by the zero set of a single polynomial. If it is not the full set $\MnR$, then it is a proper algebraic hypersurface, as we wanted to show.

We show that $M_S$ does not contain all matrices, by constructing a matrix not in it. Without loss of generality, we can assume that
\[S=\begin{pmatrix}\mathbb{I}_r&\\&-\mathbb{I}_{n-r}\end{pmatrix}\in\MnR\]
with $r<n$. Consider the matrix
\[
A=\begin{pmatrix}
0&-\mathbb{I}_r&&&\\
1&0&&&\\
&&1&&\\
&&&\ddots&\\
&&&&n-r-1
\end{pmatrix}\in\MnR
\]
whose eigenvalues are $1,\mathrm{e}^{\frac{2\pi i}{r+1}},\mathrm{e}^{\frac{4\pi i}{r+1}}\ldots, \mathrm{e}^{\frac{2r\pi i}{r+1}}$, $1,\ldots,n-r-1$. Now, for $S$ as defined above, we have that
\[
SA=\begin{pmatrix}
0&-\mathbb{I}_r&&&\\
-1&0&&&\\
&&-1&&\\
&&&\ddots&\\
&&&&-(n-r-1)
\end{pmatrix}\in\MnR.
\]
Hence the eigenvalues of $SA$ are $\mathrm{e}^{\frac{\pi i}{r+1}},\mathrm{e}^{\frac{3\pi i}{r+1}},\mathrm{e}^{\frac{5\pi i}{r+1}}\ldots, \mathrm{e}^{\frac{(2r+1)\pi i}{r+1}}$, $-1,\ldots,-(n-r-1)$, and so $A$ and $SA$ do not have common eigenvalues for a sufficiently general choice of $c$ and the $a_i$. Hence $M_S\neq \MnR$ as we wanted to show.
\end{proof}

We will be considering the map $F_{tA}$ as $t$ goes from arbitrarily small values to $1$. By~\eqref{eq:dilationaspectrum}, the map $F_{tA}$ is nondegenerate as long as $1/t\notin\specA(A)$. Let \[
1<\lambda_1<\cdots<\lambda_c
\] 
be the aligned values greater than one of $A$.
Further, let $S_1,\ldots,S_c\in\Signs$ be the corresponding sign matrices of their unique aligned trios, and $\lambda_0$ be the largest aligned value smaller than $1$ or zero if there are no such aligned values. By our initial assumption, they are not repeated. 

When $t<1/\lambda_c$, we have, by Proposition~\ref{prop:degreeFAbasics} (2), that $F_{tA}$ has degree $1$ and that $\countA(tA)=0$. We aim to show that for $t\in (1/\lambda_{k+1},1/\lambda_{k})$, we have
\[
\deg F_{tA}=1-\sum_{i=k+1}^c \sign(\chi_{AS_i}'(\lambda_i)),
\]
which gives the desired claim when $k=0$. To do so, we only have to show that the degree of $F_{tA}$ changes by $-\chi_{AS_i}'(\lambda_i)$ when $t$ passes through $1/\lambda_i$. Note that, when $t$ varies within $(1/\lambda_{k+1},1/\lambda_{k})$, neither the aligned count nor the degree changes---the latter by Proposition~\ref{prop:degreeFAbasics}(1). Hence, without loss of generality, it is enough to prove the following proposition.

\begin{proposition}\label{prop:technicalproposition}
Let $A\in\MnR$ be generic and $(1,\mathbb{I},v)$ one of its aligned trios. Assume that $1\in\specA(A)$ is a simple aligned value such that for all $S\in\mathcal{S}\setminus\{\mathbb{I}\}$, $1$ is not an eigenvalue of $SA$. Then there is an $\varepsilon>0$ such that for all $t,s\in (0,\varepsilon)$,
\[
\deg F_{(1+t)A}= \deg F_{(1-s)A} -\sign(\chi_{A}'(1)).
\]
\end{proposition}

Once this proposition is shown, we obtain the following proposition by applying the above one to $AS/\lambda$, where $(\lambda,S,v)$ is the considered aligned trio.

\begin{proposition}\label{prop:technicalproposition2}
Let $A\in\MnR$ be generic and $(\lambda,S,v)$ one of its aligned trios. Assume that $\lambda\in\specA(A)$ is a simple aligned value such that for all $T\in \in \mathcal{S}\setminus\{S\}$, $\lambda$ is not an eigenvalue of $SA$. Then there is an $\varepsilon>0$ such that for all $t,s\in (0,\varepsilon)$,
\begin{equation*}\tag*{\qed}
\deg F_{(1/\lambda+t)A}= \deg F_{(1/\lambda-s)A} -\sign(\chi_{AS}'(\lambda)).
\end{equation*}
\end{proposition}

With this proposition, the desired claim follows, i.e., that $F_{tA}$ changes the degree as wanted each time $t$ passes through the inverse of an aligned value. Note that by our original perturbation, due to Lemma \ref{lem:auxlemma}, the assumption is satisfied at each crossing.

\subsubsection*{Proof of Proposition~\ref{prop:technicalproposition}}

Since $\specSR(A)$ is discrete, there is some $\varepsilon_0>0$ such that $I_0=(1-\varepsilon_0,1+\varepsilon_0)$ does not contain other elements from $\specSR(A)$. For this interval, we have that for all $t\in I_0\setminus\{1\}$, $\mathbb{I}-tA$ is invertible; and for all $t\in I_0$ and $S\in\mathcal{S}\setminus \{\mathbb{I}\}$, $\mathbb{I}-tAS$ is invertible. If $\mathbb{I}-tAS$ is not invertible, then $1/t\in \specSR(A)$ and so, by construction of $I_0$, $t=1$ and, by assumption on the aligned value $1$, $S=\mathbb{I}$.

Now, since $1$ is an aligned value and $A$ generic, let $v\in\SnMi\cap\mathbb{R}_{>0}^n$ be the associated aligned vector with all positive entries. Now, we prove the following three lemmas in the context of Proposition~\ref{prop:technicalproposition}.

\begin{lemma}\label{lem:ex_interval}
Let $r>0$. Then there exist $x\in B_2(-v,r)$ and $\delta>0$ such that for all $t\in (1-\delta,1+\delta)\setminus\{1\}$, $x$ is a regular value of $F_{tA}$ such that
\begin{equation}\label{eq:partialdegree1}
    \sum_{\substack{z\in F_{tA}^{-1}\\z\notin \mathbb{R}^n_{\geq 0}}}\sign(\det(\partial F_{tA}(z))) 
\end{equation}
does not depend on $t$.
\end{lemma}
\begin{lemma}\label{lem:noninterseting}
There exist $r>0$ and $\delta>0$ such that for $t\in (1-\delta,1)$,
\[
B(-v,r)\cap F_{tA}\left(\bbRn_{\geq 0}\right)=\varnothing
\]
and for $t\in (1,1+\delta)$, 
\[
B(-v,r)\subset F_{tA}\left(\bbRn_{\geq 0}\right).
\]
\end{lemma}

Once these lemmas are proved, we only have to choose $r>0$, $\delta>0$ and $x\in B(-v,r)$ so that both Lemmas~\ref{lem:noninterseting} and~\ref{lem:ex_interval} apply. For this, we choose $r>0$ and $\delta>0$ as in Lemma~\ref{lem:ex_interval}, and then choose $x\in B(-v,r)$ and, if necessary, a smaller $\delta>0$, following Lemma~\ref{lem:noninterseting}. In this way, for $s\in (1-\delta,1)$, we obtain
\[
\deg F_{sA}=\sum_{\substack{z\in F_{tA}^{-1}\\z\notin \mathbb{R}^n_{\geq 0}}}\sign(\det(\partial F_{tA}(z)))
\]
and for $t\in (1,1+\delta)$,
\[
\deg F_{tA}=\sum_{\substack{z\in F_{tA}^{-1}\\z\notin \mathbb{R}^n_{\geq 0}}}\sign(\det(\partial F_{tA}(z)))+\sign(\det(\mathbb{I}-tA)),
\]
since $x$ has a preimage in the positive orthant when $t>1$. Now, $\det(\mathbb{I}-tA)=\chi_{tA}(1)$. We have that $\chi_{A}(1)=0$. Thus, for $t\in (1,1+\delta)$ sufficiently small,
\[
\sign(\det(\mathbb{I}-tA))=\sign\left(\left.\frac{d}{dt}\right|_{t=1}\chi_{tA}(1)\right).
\]
Now, $\chi_{tA}(1)=\sum_{k=0}^nc_k(A)t^{n-k}$, where $c_k(A)$ is the $k$th coefficient. Thus
\[
\left.\frac{d}{dt}\right|_{t=1}\chi_{tA}(1)=\sum_{k=0}^n(n-k)c_k(A)=n\sum_{k=0}^{n-1}c_k(A)-\sum_{k=0}^{n-1}kc_k(A)=-\sum_{k=0}^n kc_k(A)=-\chi_{A}'(1),
\]
where $\sum_{k=0}^{n-1}c_k(A)=-c_n(A)$ follows from $\chi_A(1)=0$. Hence, the desired formula follows.

We finish the proof, by proving the lemmas above. 

\begin{proof}[Proof of Lemma~\ref{lem:ex_interval}]
Denote by
\begin{equation}
    \mathcal{H}:=\{y\in\mathbb{R}^n\mid \text{there is some }i\text{ such that }x_i=0\}
\end{equation}
the union of the coordinate hyperplanes, and fix $t\in I_0\setminus \{1\}$. Then for all $x\in \mathbb{R}^n$, $x$ is a regular point of $F_{tA}$ if and only if $x\notin F_{tA}(\mathcal{H})$. If a preimage $z$ of $x$ does not lie on $\mathcal{H}$, then $F_{tA}$ is differentiable at $z$. Moreover, at that point, the Jacobian will be of the form $\mathbb{I}-tAS$, for some $S\in\Signs$, and it will be invertible since $t\in I_0\setminus \{1\}$.

Now, for all $t\in I_0$,
\[
\tilde{\mathcal{H}}_t:=\bigcup_{S\in\Signs}(\mathbb{I}-tAS)\mathcal{H}\supseteq F_{tA}(\mathcal{H})
\]
is a hyperplane arrangement. Therefore we can choose a strictly negative $x\in B_2(-v,r)$ such that $x\notin \tilde{\mathcal{H}}_1$. If we show that $I_0\ni t\mapsto \dist(x,\tilde{\mathcal{H}_t})$ is a continuous map, then, since $\dist(x,\tilde{\mathcal{H}_1})>0$, we can choose $\delta>0$ such that for all $t\in (1-\delta,1+\delta)$, $\dist(x,\tilde{\mathcal{H}}_t)>0$. Thus for all $t\in (1-\delta,1+\delta)$, $x\notin F_{tA}(\mathcal{H})$ and, by the discussion in the above paragraph, $x$ is a regular point of $F_{tA}$.

To show that $I_0\ni t\mapsto \dist(x,\tilde{\mathcal{H}_t})$ is continuous, we only have to observe that for all $i$, $S\in\Signs$ and $t\in I_0$, 
\[(\mathbb{I}-tAS)\{x\mid x_i=0\}\]
is a hyperplane, since $\{x\mid x_i\}$ does not contain any vector in the kernel of $\mathbb{I}-tAS$. Only $\mathbb{I}-A$ has a kernel, but it is a line spanned by the strictly positive vector $v$. Moreover, taking the generalized cross product of the $(\mathbb{I}-tAS)e_j$, $j\neq i$, and normalizing it, we can obtain a continuous map
\[I_0\ni t\mapsto n(S,i,t)\in\mathbb{S}^{n-1}\]
that maps $t$ to a normal vector of $(\mathbb{I}-tAS)\{x\mid x_i=0\}$. Hence
\[
\dist(x,\tilde{\mathcal{H}}_t)=\min\left\{|n(S,i,t)^tx|\mid i\in\{1,\ldots,n\},\,S\in\Signs\right\}
\]
is a continuous function as we wanted to show.

Now, for all $t\in (1-\delta,1+\delta)$ and all $S\in\Signs\setminus \{1\}$, we have that
\[
(\mathbb{I}-tAS)^{-1}x\notin \mathcal{H}.
\]
Note that, by choosing $\delta>0$ small enough if necessary, we can guarantee $t\in I_0$ and that $\mathbb{I}-tAS$ is invertible. Hence we have that the signs of $(\mathbb{I}-tAS)^{-1}x$ are constant. This means that, independently of $t\in (1-\delta,1+\delta)$,
\[
F_{tA}^{-1}(x)\cap S\mathbb{R}_{> 0}^n
\]
is either empty or has size one. Now, for $z\in S\mathbb{R}_{> 0}^n$, $\sign(\det(\delta F_{tA}(z)))=\sign(\det(\mathbb{I}-tAS))$ is constant. Hence the sum \eqref{eq:partialdegree1} remains constant as desired.
\end{proof}
For the proof of Lemma~\ref{lem:noninterseting}, we will need the following auxiliary lemma.
\begin{lemma}\label{lem:trivialintersection}
Let $A\in\MnR$. If $\lambda$ is a simple eigenvalue of $A$, then
\[
	\ker(\lambda\mathbb{I}-A) \cap \img (\lambda\mathbb{I}-A)  =0.
\]
\end{lemma}
\begin{proof}[Proof of Lemma~\ref{lem:trivialintersection}]
If $v\in\ker(\lambda\mathbb{I}-A)\cap \img(\lambda\mathbb{I}-A)$ is non-zero, then, by taking $w$ such that $v=(\lambda\mathbb{I}-A)w$, there is $w\neq 0$ such that $w\in \ker(\lambda\mathbb{I}-A)^2$, but $w\notin \ker(\lambda\mathbb{I}-A)$. Thus $w$ is a nontrivial generalized eigenvector of rank 2 of $A$ corresponding to $\lambda$. However, $\lambda$ is a simple eigenvalue of $A$, so this is impossible.
\end{proof}

\begin{proof}[Proof of Lemma \ref{lem:noninterseting}]
By Lemma \ref{lem:trivialintersection}, we have that
\[-v\notin (\mathbb{I}-A)\mathcal{H},\]
where $\mathcal{H}$ is the union of the coordinate hyperplanes, as in the proof of Lemma~\ref{lem:ex_interval}. Arguing as in that proof, we have that $I_0\ni t\mapsto \dist(-v,(\mathbb{I}-tA)\mathcal{H})$ is continuous, and so, since $\dist(-v,(\mathbb{I}-A)\mathcal{H})$ is positive, we have that there is $\delta>0$ and $r>0$ such that for all $t\in(1-\delta,1+\delta)\subseteq I_0$,
\[\dist(-v,(\mathbb{I}-tA)\mathcal{H})>r.\]
We show now that these are the desired $r$ and $\delta$.

Fix $t\in (1-\delta,1)$. $F_{tA}(\bbRn_{\geq 0})$ is a closed pointed cone, since it is image of a closed pointed cone under the invertible linear map $\mathbb{I}-tA$. This cone contains $v$, so it does not contain $-v$. Now, 
\[
\dist(-v,F_{tA}(\bbRn_{\geq 0}))=\dist(-v,\mathrm{bd}F_{tA}(\bbRn_{\geq 0}))\geq \dist(-v,(\mathbb{I}-A)\mathcal{H}),
\]
since the nearest point in $F_{tA}(\bbRn_{\geq 0})$ to $-v$ lies in the boundary, $\mathrm{bd}F_{tA}(\bbRn_{\geq 0})$ of $F_{tA}(\bbRn_{\geq 0})$, which is contained inside $(\mathbb{I}-A)\mathcal{H}$. Hence, by the first paragraph in the proof, $\dist(-v,F_{tA}(\bbRn_{\geq 0}))>r$, and so $B(-v,r)\cap F_{tA}(\bbRn_{\geq 0}))=\varnothing$, as desired.

Fix $t\in (1,1+\delta)$. $F_{tA}(\bbRn_{>0})$ is a full-dimensional cone, since it is the image of a full-dimensional cone under the invertible linear map $\mathbb{I}-tA$. Now, we have that
\[
\dist(-v,\mathrm{bd}F_{tA}(\bbRn_{\geq 0}))\geq \dist(-v,(\mathbb{I}-A)\mathcal{H}).
\]
Therefore, by the first paragraph in the proof, $\dist(-v,\mathrm{bd}F_{tA}(\bbRn_{\geq 0}))>r$, and so $B(-v,r)\subseteq F_{tA}(\bbRn_{\geq 0}))$, as desired.
\end{proof}

\printbibliography
\end{document}